\documentclass[12pt]{article}

\usepackage[authoryear]{natbib}
\usepackage{graphicx}
\usepackage{amssymb}
\usepackage{amsmath}
\usepackage{times}
\usepackage{bm}

\newcommand{\Es}{\mathbb{E}}
\newcommand{\diag}{{\rm diag}}
\newcommand{\dd}{{\rm d}}
\newcommand{\tr}{\textrm{tr}}

\title{Local power of the LR, Wald, score and gradient tests in dispersion models}
\author{Artur J.~Lemonte,\quad Silvia L.P.~Ferrari\\
{\small {\em Departamento de Estat\'istica, Universidade de S\~ao Paulo, Brazil}}}
\date{}

\begin{document}
\maketitle

\begin{abstract}

We derive asymptotic expansions up to
order $n^{-1/2}$ for the nonnull distribution functions
of the likelihood ratio, Wald, score and gradient
test statistics in the class of dispersion models,
under a sequence of Pitman alternatives.
The asymptotic distributions of these statistics are
obtained for testing a subset of regression
parameters and for testing the precision parameter.
Based on these nonnull asymptotic expansions it is shown that
there is no uniform superiority of one test
with respect to the others for testing a subset of regression
parameters. Furthermore, in order
to compare the finite-sample performance of these tests
in this class of models, Monte Carlo simulations are presented. An empirical
application to a real data set is considered for
illustrative purposes.\\

\end{abstract}

\section{Introduction}\label{introduction}

The paper by \cite{NelderWedderburn1972} introduced
the class of generalised linear models (GLMs) and showed that a large
variety of non-normal data may be analysed by a simple general
technique \citep[see, for example,][]{McCu:Neld:1989,
DobsonBarnett2008}. The GLMs were originally developed for the
exponential family of distributions, but the main ideas
were extended to a wider class of models called dispersion models (DMs)
in such a way that most of their good properties were preserved.
This class of models was introduced by \cite{Jorgensen1987a} and
studied in details in \cite{Jorgensen1997a}. Some recent references
about DMs are \cite{Kokonendji-et-al-2004},
\cite{Jorgensen-et-al-2010}, \cite{Simas-et-al-2010} and
\cite{Rocha-et-al-2010}.

The class of DMs with position parameter $\theta$ (which vary in an interval
of the real line) and precision parameter $\phi>0$ has
probability density function of the form
\begin{equation}\label{dens1}
\pi(y;\theta,\phi) = \exp\{\phi t(y,\theta) + c(y,\phi)\}, 
\end{equation}
where $t(\cdot,\cdot)$ and $c(\cdot,\cdot)$ are known functions.
If $Y$ is continuous, $\pi$ is assumed to be a density with
respect to the Lebesgue measure, while if $Y$ is discrete, $\pi$
is assumed to be a density with respect to the counting measure.
The parameter $\theta$ may be generally interpreted as a kind of location parameter,
not necessarily the mean of the distribution. Several
models of the form (\ref{dens1}) are discussed by \cite{Jorgensen1987a,
Jorgensen1987b, Jorgensen1997a}, who also examined
their statistical properties. It is evident that
some special cases arise from (\ref{dens1}). Exponential dispersion models (EDMs)
represent a special case of DMs with $t(y,\theta)=y\theta-b(\theta)$,
where $\Es(Y)=\dd b(\theta)/\dd\theta$; see \cite{Jorgensen1992}.
An important subclass of DMs of
special interest, called proper dispersion models (PDMs),
arise when $c(y,\phi)$ is additive,
i.e.~$c(y,\phi)=a_1(y)+a_2(\phi)$, where $a_1(\cdot)$
and $a_2(\cdot)$ are known functions \citep[see, for instance,][]{Jorgensen1997b}.
The class of PDMs covers important distributions which are not covered by the EDMs,
such as the log-gamma distribution, the McCullagh distribution
\citep{McCullagh1989}, the reciprocal inverse Gaussian distribution
and the simplex distribution, which is suitable for modeling continuous proportions
\citep{BNJorgensen1991}. The von
Mises distribution, which also belongs to the class of PDMs and does
not belong to the EDMs, is particularly useful for the analysis
of circular data; see \cite{MardiaJupp2000}. The PDMs have two important
general properties. First, the distribution of the statistic $T = t(Y,\theta)$
does not depend on $\theta$ when $\phi$ is known, that is, $T$
is a pivotal quantity for $\theta$. Second, (\ref{dens1}) is an
exponential family with canonical statistic $T$ when $\theta$
is known.

Large-sample tests, such as the likelihood ratio, Wald and Rao score tests,
are usually employed for testing hypotheses in parametric models.
A new criterion for testing hypotheses, referred to as the  {\it gradient test},
was proposed in \cite{Terrell2002}.
Its statistic is very simple to compute when compared with
the other three classic statistics. Here, it is worthwhile to quote 
\cite{Rao2005}: ``The suggestion by Terrell is attractive as it is simple to compute.
It would be of interest to investigate the performance of the [gradient] statistic.''
Also, Terrell's statistic shares the same first
order asymptotic properties with the
likelihood ratio, Wald and score statistics.
That is, to the first order of approximation, the likelihood ratio,
Wald, score and gradient statistics have the same asymptotic
distributional properties either under the null hypothesis or under a sequence of
Pitman alternatives, i.e.~a sequence of local
alternatives that shrink to the null hypothesis at a convergence rate $n^{-1/2}$.
Additionally, it is known that, up to an error of order $n^{-1}$, the likelihood ratio,
Wald, score and gradient tests have the
same size properties but their local powers differ in the $n^{-1/2}$ term.
Therefore, a meaningful comparison among the criteria can be performed by comparing
the nonnull asymptotic expansions to order $n^{-1/2}$ of the distribution functions
of these statistics under a sequence of Pitman alternatives.

In this paper, our main objective is to derive
nonnull asymptotic expansions to order $n^{-1/2}$
of the distribution functions of the likelihood ratio, Wald, score and
gradient statistics under a sequence of local alternatives
and to compare the local power of the corresponding tests in the class of DMs.
In order to compare the finite-sample performance of these tests in this class of models we
also perform a Monte Carlo simulation study.
As far as we know, there is no mention in the
statistical literature on the use of the gradient test in DMs.

The nonnull asymptotic expansions up to order $n^{-1/2}$
for the distribution functions of the likelihood ratio and Wald statistics were derived by
\cite{Hayakawa1975}, while an analogous result for the
score statistic was obtained by \cite{HarrisPeers1980}.
The asymptotic expansion up
to order $n^{-1/2}$ for the distribution functions of the gradient statistic was derived
by \cite{LemFer2010}. The expansions are
very general, although being difficult or even impossible to
particularize their formulas for specific regression models.
As we shall see below, we have been capable to apply their results for DMs.
In particular, we derive closed-form expressions for the
coefficients that define the nonnull asymptotic expansions
of these statistics in this class of models and show
that there is no uniform superiority
of one test with respect to the others for testing a subset of
regression parameters.

The rest of the paper is organized as follows. Section \ref{tests} briefly describes the
likelihood ratio, Wald, score and gradient tests. We present the class
of DMs in Section \ref{DMs}. In Section \ref{main_result} we derive
the nonnull asymptotic expansions of the likelihood ratio,
Wald, score and gradient statistics for testing a subset
of regression parameters in DMs. The local power of the
likelihood ratio, Wald, score and gradient tests are compared in
Section \ref{power_comparison}. In Section \ref{test_phi}
we consider hypothesis testing on the precision parameter.
Monte Carlo simulation results are addressed in Section \ref{simulations}.
We consider an empirical application in Section \ref{applications}
for illustrative purposes. Section \ref{conclusions} closes the
paper with some concluding remarks.

\section{Background}\label{tests}

Let $\ell(\bm{\theta})$, $\bm{U}_{\bm{\theta}}$ and $\bm{K}_{\bm{\theta}}$
denote the total log-likelihood function, the score
function and the information matrix for the parameter vector
$\bm{\theta}=(\theta_{1},\ldots,\theta_{k})^\top$ of dimension $k$, respectively.
Let $\bm{K}_{\bm{\theta}}^{-1}$ denote the inverse of $\bm{K}_{\bm{\theta}}$.
Consider the partition $\bm{\theta} = (\bm{\theta}_{1}^{\top}, \bm{\theta}_{2}^{\top})^{\top}$,
where the dimensions of $\bm{\theta}_{1}$ and $\bm{\theta}_{2}$
are $q$ and $k-q$, respectively.
Suppose the interest lies in testing the composite null hypothesis
$\mathcal{H}_{0}:\bm{\theta}_{2} = \bm{\theta}_{20}$
against $\mathcal{H}_{1}:\bm{\theta}_{2}\neq\bm{\theta}_{20}$,
where $\bm{\theta}_{20}$ is a specified vector.
Hence, $\bm{\theta}_{1}$ acts as a vector of nuisance parameters.
The likelihood ratio ($S_1$), Wald ($S_2$), score ($S_3$) and gradient ($S_4$) statistics for
testing $\mathcal{H}_{0}$ versus $\mathcal{H}_{1}$ are given, respectively, by
\[
S_{1} = 2\bigl\{\ell(\widehat{\bm{\theta}}) - \ell(\widetilde{\bm{\theta}})\bigr\},
\qquad
S_{2} = (\widehat{\bm{\theta}} - \widetilde{\bm{\theta}})^{\top}\widehat{\bm{K}}_{\bm{\theta}}
(\widehat{\bm{\theta}} - \widetilde{\bm{\theta}}),
\]
\[
S_{3} = \widetilde{\bm{U}}_{\bm{\theta}}^\top\widetilde{\bm{K}}_{\bm{\theta}}^{-1}\widetilde{\bm{U}}_{\bm{\theta}},
\qquad
S_{4} = \widetilde{\bm{U}}_{\bm{\theta}}^\top(\widehat{\bm{\theta}} - \widetilde{\bm{\theta}}),
\]
where $\widehat{\bm{\theta}}=(\widehat{\bm{\theta}}_1^\top,\widehat{\bm{\theta}}_2^\top)^\top$ and
$\widetilde{\bm{\theta}}=(\widetilde{\bm{\theta}}_1^\top,\bm{\theta}_{20}^{\top})^\top$
denote the maximum likelihood estimators of
$\bm{\theta} =(\bm{\theta}_{1}^{\top}, \bm{\theta}_{2}^{\top})^{\top}$ under
$\mathcal{H}_{1}$ and $\mathcal{H}_{0}$, respectively,
$\widehat{\bm{K}}_{\bm{\theta}}=\bm{K}_{\bm{\theta}}(\widehat{\bm{\theta}})$,
$\widetilde{\bm{K}}_{\bm{\theta}}=\bm{K}_{\bm{\theta}}(\widetilde{\bm{\theta}})$ and
$\widetilde{\bm{U}}_{\bm{\theta}} = \bm{U}_{\bm{\theta}}(\widetilde{\bm{\theta}})$.
The limiting distribution of $S_1$, $S_2$, $S_3$ and $S_4$ is $\chi_{k-q}^2$ under
$\mathcal{H}_{0}$ and $\chi_{k-q,\lambda}^{2}$, i.e.~a noncentral chi-square
distribution with $k-q$ degrees of freedom and an appropriate noncentrality parameter $\lambda$,
under $\mathcal{H}_{1}$. The null hypothesis is rejected for a given nominal level, $\gamma$ say,
if the test statistic exceeds the upper $100(1-\gamma)\%$ quantile of the $\chi_{k-q}^{2}$
distribution.

From the partition of $\bm{\theta}$, we have the corresponding partitions
\[
\bm{U}_{\bm{\theta}} = (\bm{U}_{\bm{\theta}_1}^\top, \bm{U}_{\bm{\theta}_2}^\top)^\top,
\quad
\bm{K}_{\bm{\theta}} =
\begin{bmatrix}
\bm{K}_{11} & \bm{K}_{12} \\
\bm{K}_{21} & \bm{K}_{22}
\end{bmatrix},
\quad
\bm{K}_{\bm{\theta}}^{-1} =
\begin{bmatrix}
\bm{K}^{11} & \bm{K}^{12} \\
\bm{K}^{21} & \bm{K}^{22}
\end{bmatrix}.
\]
Thus, the statistics $S_2$, $S_3$ and $S_4$ can be rewritten as
\[
S_{2} = (\widehat{\bm{\theta}}_{2} - \bm{\theta}_{20})^{\top}\widehat{\bm{K}}^{22^{-1}}
(\widehat{\bm{\theta}}_{2} - \bm{\theta}_{20}),
\quad
S_{3} = \widetilde{\bm{U}}_{\bm{\theta}_2}^\top\widetilde{\bm{K}}^{22}\widetilde{\bm{U}}_{\bm{\theta}_2},
\quad
S_{4} = \widetilde{\bm{U}}_{\bm{\theta}_2}^\top(\widehat{\bm{\theta}}_{2} - \bm{\theta}_{20}),
\]
where $\widehat{\bm{K}}^{22}=\bm{K}^{22}(\widehat{\bm{\theta}})$,
$\widetilde{\bm{K}}^{22}=\bm{K}^{22}(\widetilde{\bm{\theta}})$ and
$\widetilde{\bm{U}}_{\bm{\theta}_2} = \bm{U}_{\bm{\theta}_2}(\widetilde{\bm{\theta}})$.

Noticed that $S_{4}$ has a very simple form and does not involve
the information matrix, neither expected nor observed, unlike $S_{2}$ and $S_{3}$.
\cite{Terrell2002} points out that the
gradient statistic ``is not transparently non-negative, even
though it must be so asymptotically.'' His Theorem 2 implies that
if the log-likelihood function is concave and is differentiable at
$\widetilde{\bm{\theta}}$, then $S_{4}\ge 0$.

Recently, \cite{LemonteFerrari2011}
obtained the nonnull asymptotic expansions of the
likelihood ratio, Wald, score and gradient statistics in
Birnbaum--Saunders regression models \citep{RiekNedelman91}.
An interesting finding is that, up to an error of order $n^{-1}$,
the four tests have the same local power in this class of models.
Their simulation study evidenced that the score and the gradient tests
perform better than the likelihood ratio and Wald tests in small and
moderate-sized samples and hence they concluded that the
gradient test is an appealing alternative to the three
classic asymptotic tests in Birnbaum--Saunders regressions.

\section{Dispersion models}\label{DMs}

We assume that the random variables $y_1,\ldots, y_n$ are independent
and each $y_l$ has a probability density function of the form
\begin{equation}\label{dens2}
\pi(y_{l};\theta_{l},\phi) = \exp\{\phi t(y_{l},\theta_{l}) + c(y_{l},\phi)\},
\quad l = 1,\ldots,n.
\end{equation}
The mean of $Y_l$ will be denoted by $\mu_{l}$, and is not
necessary equal to $\theta_{l}$, the parameter of interest.
In order to introduce a regression structure in the class of models in (\ref{dens2}),
we assume that
\begin{equation}\label{sistpart}
d(\theta_{l}) = \eta_{l} = f(\bm{x}_l;\bm{\beta}),\quad l = 1,\ldots,n,
\end{equation}
where $d(\cdot)$ is a known one-to-one differentiable link function,
$\bm{x}_{l}=(x_{l1},\ldots,x_{lm})^\top$ is an $m$-vector of nonstocastic
variables associated with the $l$-th response, $\bm{\beta}=(\beta_1,\ldots,\beta_p)^\top$
is a set of unknown parameters to be estimated ($m\leq p<n$), and $f(\cdot;\cdot)$ is a possible
nonlinear twice continuous differenciable function with respect to $\bm{\beta}$.
The regression structure links the covariates $\bm{x}_l$ to the parameter
of interest $\theta_{l}$. The
$n\times p$ matrix of derivatives of $\bm{\eta}=(\eta_1,\ldots,\eta_n)^\top$ with
respect to $\bm{\beta}$, specified by $\bm{X}^*=\partial\bm{\eta}/\partial\bm{\beta}^\top$,
is assumed to be of full rank, i.e.~rank$(\bm{X}^*)=p$ for all $\bm{\beta}$.
Further, it is assumed that the precision parameter is
unknown and it is the same for all observations. It is also assumed that the
usual regularity conditions for maximum likelihood estimation and large sample
inference hold; see \citet[Ch.~9]{CoxHinkley1974}.

The class of regression models defined by
(\ref{dens2}) and (\ref{sistpart}) extends the class of generalised
linear models discussed by \cite{McCu:Neld:1989} in two directions.
First and as noted before, it includes important distributions
which are not exponential family models. Second, it allows for a nonlinear
structure in $\bm{\eta}$. The class of models in
(\ref{dens2})-(\ref{sistpart}) is also a natural extension
of the exponential family nonlinear models (EFNLMs)
introduced by  \cite{CordPaula1989}, which in turn
extends the well-known GLMs by allowing the regression
structure to be nonlinear. The EFNLMs are defined by
equations (\ref{dens2}) and (\ref{sistpart}), with $t(y_{l},\theta_{l})=y_{l}\theta_{l}-b(\theta_{l})$
and $c(y_{l},\phi)=a_1(y_{l})+a_2(\phi)$ in (\ref{dens2}).

Let $\ell=\ell(\bm{\beta},\phi)=\sum_{l=1}^{n}\{\phi t(y_{l},\theta_{l}) + c(y_{l},\phi)\}$
be the total log-likelihood function
for $\bm{\beta}$ and $\phi$, where $\theta_{l}$ is related to $\bm{\beta}$ by (\ref{sistpart}). We define
$D_{il}=D_{il}(\theta_{l},\phi)=\Es\{\partial^{i}t(Y_l,\phi)/\partial\theta_{l}^i\}$,
for $i=1,2,3$ and $l=1,\ldots,n$. From regularity conditions we have that $D_{1l}=0$,
for $l=1,\ldots,n$. Table \ref{tab1}
lists $D_{2l}$ and $D_{3l}$ for some dispersion models.
The total score function and the total Fisher information matrix for $\bm{\beta}$
are given, respectively, by
$\bm{U}_{\bm{\beta}}=\phi\bm{X}^{*\top}\dot{\bm{t}}$ and
$\bm{K}_{\bm{\beta}} = \phi\bm{X}^{*\top}\bm{W}\bm{X}^*$,
where $\dot{\bm{t}}=\dot{\bm{t}}(\bm{y},\bm{\theta})=(\dot{t}_1,\ldots,\dot{t}_n)^\top$ 
is an $n\times 1$ vector with $\dot{t}_l=\partial t(y_l,\theta_{l})/\partial\theta_{l}$,
$\bm{y}=(y_1,\ldots,y_n)^\top$,
$\bm{\theta}=(\theta_{1},\ldots,\theta_{n})^\top$ and
$\bm{W}=\diag\{w_{1},\ldots,w_{n}\}$ with $w_{l} = -D_{2l}(\dd\theta_{l}/\dd\eta_{l})^2$.
A simple calculation shows that
$\Es(\partial^2\ell/\partial\bm{\beta}\partial\phi)=\bm{0}$ and then the parameters
$\bm{\beta}$ and $\phi$ are globally orthogonal \citep{CoxReid1987}.
Let $\alpha_i=\sum_{l=1}^{n}\Es\{\partial^ic(Y_l,\phi)/\partial\phi^i\}
=\sum_{l=1}^{n}\Es\{c^{(i)}(Y_l,\phi)\}$, for $i=1,2,3$. The derivatives
of the $\alpha_i$'s with respect to $\phi$ are written with primes,
i.e.~$\alpha_{i}' = \dd\alpha_i/\dd\phi$ and so on. We have that
the joint information matrix for $(\bm{\beta}^\top,\phi)^\top$ is given by
$\diag\{\bm{K}_{\bm{\beta}},-\alpha_2\}$.
\begin{table}[!htp]
{\footnotesize
\begin{center}
\caption{Expressions of $D_{2l}$ and $D_{3l}$ ($l=1,\ldots,n$) for some dispersion models.$^{\dagger}$}\label{tab1}
\begin{tabular}{lcc}\hline
Model            &  $D_{2l}$                                   & $D_{3l}$     \\\hline
Normal           &  $-1$                          & 0    \\
Inverse Gaussian &  $-(-2\theta_{l})^{-3/2}$                          & $-3(-2\theta_{l})^{-5/2}$   \\
Reciprocal inverse Gaussian &   $-1/\theta_{l}$                          & 0   \\
Gamma            &  $-1/\theta_{l}^2$ & $2/\theta_{l}^3$ \\
Reciprocal gamma &  $-1/\theta_{l}^2$ & $2/\theta_{l}^3$ \\
Log-gamma        &  $-1$ & 1 \\
von Mises        &  $-I_1(\phi)/I_0(\phi)$ & 0\\
generalised hyperbolic secant & $2/(\theta_l^2 + 1)^3$ & $(2\theta_l^3 + 10\theta_l)/(\theta_l^2 + 1)^3$\\\hline
\multicolumn{3}{l}{{\footnotesize $^{\dagger}$$I_j(\phi)$ is the modified Bessel function of the first kind and order $j$.}}
\end{tabular}
\end{center}    }
\end{table}

The maximum likelihood estimate
(MLE) $\widehat{\bm{\beta}}$ of $\bm{\beta}$ can be obtained iteratively using
standard reweighted least squares method \citep{Jorgensen1983,Jorgensen1984}:
\[
\bm{X}^{*(m)\top}\bm{W}^{(m)}\bm{X}^{*(m)}\bm{\beta}^{(m+1)} = \bm{X}^{*(m)\top}\bm{W}^{(m)}\bm{y}^{*(m)},
\quad m = 0,1,\ldots,
\]
where $\bm{y}^{*(m)} = \bm{X}^{*(m)}\bm{\beta}^{(m)} + \bm{N}^{(m)}\dot{\bm{t}}^{(m)}$
is an adjusted dependent variable and $\bm{N}$ is a diagonal matrix
given by $\bm{N}=-\diag\{D_{21}^{-1}(\dd\theta_{1}/\dd\eta_{1})^{-1}, \ldots,
D_{2n}^{-1}(\dd\theta_{n}/\dd\eta_{n})^{-1}\}$. The estimate $\widehat{\bm{\beta}}$ depends
directly on the distribution only through the function $D_{2l}$ and
does not depend on the parameter $\phi$. The maximum likelihood
estimate $\widehat{\phi}$ of $\phi$ is the solution of
\begin{equation}\label{est_phi}
\sum_{l=1}^{n}\{t(y_l,\widehat{\theta}_l) + c^{(1)}(y_l,\widehat{\phi})\} = 0.
\end{equation}
The maximum likelihood estimators $\widehat{\bm{\beta}}$
and $\widehat{\phi}$ are asymptotically independent due to
their asymptotic normality and the block diagonal structure of the joint
information matrix. If the model is a PDM the $\alpha_{i}$'s
can be expressed as functions of $\phi$ only, namely $\alpha_{i}=na_2^{(i)}(\phi)$
for $i=1,2,3$, where $a_2^{(i)}(\phi)$ is the $i$-th derivative of $a_2(\phi)$ with
respect to $\phi$. In this case, the $(p+1,p+1)$-th element of the joint
information matrix is simply $-na_2^{(2)}(\phi)$ and equation
(\ref{est_phi}) reduces to $a_2^{(1)}(\widehat{\phi})=-\sum_{l=1}^{n}t(y_l,\widehat{\theta}_l)/n$.

In what follows, we shall consider tests based on the
likelihood ratio ($S_1$), Wald ($S_2$), Rao score ($S_3$) and gradient ($S_4$) statistics
in the class of DMs for testing a composite null hypothesis
$\mathcal{H}_{0}: \bm{\beta}_{2} = \bm{\beta}_{20}$.
This hypothesis will be tested against the alternative hypothesis
$\mathcal{H}_{1}:\bm{\beta}_{2}\neq\bm{\beta}_{20}$,
where $\bm{\beta}$ is partitioned as $\bm{\beta} = (\bm{\beta}_{1}^{\top},
\bm{\beta}_{2}^{\top})^{\top}$, with
$\bm{\beta}_{1} = (\beta_{1}, \dots,\beta_{q})^{\top}$ and
$\bm{\beta}_{2} = (\beta_{q+1},\dots,\beta_{p})^{\top}$. Here,
$\bm{\beta}_{20}$ is a fixed column vector of dimension $p-q$. The
partition of the parameter vector $\bm{\beta}$ induces the corresponding partitions
$\bm{U}_{\bm{\beta}} = (\bm{U}_{\bm{\beta}_1}^\top, \bm{U}_{\bm{\beta}_2}^\top)^\top$,
with $\bm{U}_{\bm{\beta}_1}=\phi\bm{X}_1^{*\top}\dot{\bm{t}}$ and
$\bm{U}_{\bm{\beta}_2}=\phi\bm{X}_2^{*\top}\dot{\bm{t}}$,
\[
\bm{K}_{\bm{\beta}} =
\begin{bmatrix}
\bm{K}_{\bm{\beta}11} & \bm{K}_{\bm{\beta}12} \\
\bm{K}_{\bm{\beta}21} & \bm{K}_{\bm{\beta}22}
\end{bmatrix} = \phi
\begin{bmatrix}
\bm{X}_1^{*\top}\bm{W}\bm{X}_1^{*} & \bm{X}_1^{*\top}\bm{W}\bm{X}_2^{*} \\
\bm{X}_2^{*\top}\bm{W}\bm{X}_1^{*} & \bm{X}_2^{*\top}\bm{W}\bm{X}_2^{*}
\end{bmatrix},
\]
with the matrix $\bm{X}^*$ partitioned
as $\bm{X}^* = \bigl[\bm{X}_{1}^*\ \ \bm{X}_{2}^*\bigr]$,
$\bm{X}_{1}^*$ being $n\times q$ and $\bm{X}_{2}^*$ being $n\times (p-q)$.
Let $(\widehat{\bm{\beta}}_{1}, \widehat{\bm{\beta}}_{2}, \widehat{\phi})$ and
$(\widehat{\bm{\beta}}_{1}, \bm{\beta}_{20}, \widetilde{\phi})$ be the
unrestricted and restricted MLEs of $(\bm{\beta}_{1}, \bm{\beta}_{2}, \phi)$, respectively.
The likelihood ratio, Wald, score and gradient statistics for testing $\mathcal{H}_{0}$
can be expressed, respectively, as
\[
S_{1} = 2\bigl\{\ell(\widehat{\bm{\beta}}_{1}, \widehat{\bm{\beta}}_{2},\widehat{\phi})
- \ell(\widetilde{\bm{\beta}}_{1}, \bm{\beta}_{20},\widetilde{\phi})\bigr\},
\qquad
S_{2} = \widehat{\phi}(\widehat{\bm{\beta}}_{2}- \bm{\beta}_{20})^{\top}
(\widehat{\bm{R}}^{\top}\widehat{\bm{W}}\widehat{\bm{R}})
(\widehat{\bm{\beta}}_{2} - \bm{\beta}_{20}),
\]
\[
S_{3} = \widetilde{\bm{s}}^{\top}\widetilde{\bm{W}}^{1/2}\widetilde{\bm{X}}_{2}^*
(\widetilde{\bm{R}}^{\top}\widetilde{\bm{W}}\widetilde{\bm{R}})^{-1}
\widetilde{\bm{X}}_{2}^{*\top}\widetilde{\bm{W}}^{1/2}\widetilde{\bm{s}},\qquad
S_{4} = \widetilde{\phi}^{1/2}\widetilde{\bm{s}}^{\top}\widetilde{\bm{W}}^{1/2}\widetilde{\bm{X}}_{2}^*
(\widehat{\bm{\beta}}_{2} - \bm{\beta}_{20}),
\]
where $\bm{s} = (s_1,\ldots,s_n)^\top$ with $s_l=\phi^{1/2}\dot{t}_{l}(-D_{2l})^{-1/2}$ and
$\bm{R} = \bm{X}_{2}^* - \bm{X}_{1}^*(\bm{X}_{1}^{*\top}\bm{W}\bm{X}_{1}^*)^{-1}
\bm{X}_{1}^{*\top}\bm{W}\bm{X}_{2}^*$.
Here, tildes and hats indicate evaluation at the restricted and unrestricted MLEs,
respectively. The limiting distribution of all these statistics
under $\mathcal{H}_{0}$ is $\chi_{p-q}^2$. Note that,
unlike the Wald and score statistics, the gradient statistic
does not involve any matrix inversion.

\section{Nonnull asymptotic distributions in DMs}\label{main_result}

We present in this section expressions for the nonnull asymptotic
expansions up to order $n^{-1/2}$ for the nonnull distribution of
the likelihood ratio, Wald, score and gradient statistics for
testing a subset of regression parameters in DMs.
It should be mentioned that the general nonnull asymptotic expansions derived in
\cite{Hayakawa1975}, \cite{HarrisPeers1980} and \cite{LemFer2010}
were developed for continuous distributions. It implies that the results derived in this
section are only valid for continuous DMs. Here, we shall assume the following local alternative
hypothesis $\mathcal{H}_{1n}:\bm{\beta}_{2}=\bm{\beta}_{20} + \bm{\epsilon}$,
where $\bm{\epsilon} = (\epsilon_{q+1}, \ldots,\epsilon_{p})^{\top}$ with
$\epsilon_{r} = O(n^{-1/2})$ for $r=q+1,\ldots,p$.

We introduce the following quantities:
\[
\bm{\epsilon}^{*} = \begin{bmatrix}
\bm{K}_{\bm{\beta}11}^{-1}\bm{K}_{\bm{\beta}12}\\
-\bm{I}_{p-q}
\end{bmatrix}\bm{\epsilon},
\quad
\bm{A} = \begin{bmatrix}
                       \bm{K}_{\bm{\beta}11}^{-1} & \bm{0}\\
                       \bm{0} & \bm{0}
                      \end{bmatrix},
\quad
\bm{M} = \bm{K}_{\bm{\beta}}^{-1} - \bm{A},
\]
where $\bm{I}_{p-q}$ is a $(p-q)\times(p-q)$ identity matrix.
Additionally, let $\bm{Z} = \bm{X}^*(\bm{X}^{*\top}\bm{W}\bm{X}^*)^{-1}\bm{X}^{*\top} = \{z_{lm}\}$,
$\bm{Z}_{1} = \bm{X}_{1}^*(\bm{X}_{1}^{*\top}\bm{W}\bm{X}_{1}^*)^{-1}\bm{X}_{1}^{*\top} = \{z_{1lm}\}$,
\[
\bm{X}_{l}^* = \biggl\{\frac{\partial^2\eta_{l}}{\partial\beta_{r}\partial\beta_{s}}\biggr\} =
\begin{bmatrix}
\bm{X}_{11l}^* & \bm{X}_{12l}^* \\
\bm{X}_{21l}^* & \bm{X}_{22l}^*
\end{bmatrix},\quad r,s=1,\ldots,p,\quad l =1,\ldots,n,
\]
$\bm{Z}_{d} = \diag\{z_{11},\ldots, z_{nn}\}$,
$\bm{Z}_{1d} = \diag\{z_{111},\ldots, z_{1nn}\}$,
$\bm{F} = \diag\{f_{1},\ldots,f_{n}\}$, $\bm{G} = \diag\{g_{1},\ldots,g_{n}\}$,
$\bm{E} = \diag\{e_{1},\ldots,e_{n}\}$,
$\bm{t} = (t_{1},\ldots,t_{n})^{\top} = \bm{X}^*\bm{\epsilon}^{*}$,
$\bm{b} = (b_{1},\ldots,b_{n})^{\top} = \bm{X}_{2}^*\bm{\epsilon}$,
$\bm{T} = \diag\{t_{1},\ldots,t_{n}\}$,
$\bm{T}^{(2)} = \bm{T}\odot\bm{T}$, $\bm{T}^{(3)} = \bm{T}^{(2)}\odot\bm{T}$
and $\bm{B} = \diag\{b_{1},\ldots,b_{n}\}$, where ``$\odot$'' denotes the Hadamard (direct)
product of matrices, and
\[
f_l = -\frac{\dd \theta_l}{\dd \eta_l}\frac{\dd ^2\theta_l}{\dd \eta_l^2}D_{2l}
            - \biggl(\frac{\dd\theta_l}{\dd \eta_l}\biggr)^{3}D_{3l},\quad
g_l = -\frac{\dd\theta_l}{\dd \eta_l}\frac{\dd ^2\theta_l}{\dd\eta_l^2}D_{2l},\quad
e_l = - \biggl(\frac{\dd\theta_l}{\dd \eta_l}\biggr)^{3}D_{2l}',
\quad l=1,\ldots,n,
\]
where $D_{2l}'$ denotes the first derivative of
$D_{2l}$ with respect to $\theta_{l}$, for $l=1,\ldots,n$.

The nonnull distributions of $S_1$, $S_2$, $S_3$ and $S_4$
under Pitman alternatives for testing $\mathcal{H}_{0}:\bm{\beta}_{2}=\bm{\beta}_{20}$
in DMs can be expressed as
\begin{equation*}\label{expansion}
\Pr(S_{i}\leq x) = G_{p-q,\lambda}(x) + \sum_{k=0}^{3}b_{ik}G_{p-q+2k,\lambda}(x) + O(n^{-1}),
\quad i=1,2,3,4,
\end{equation*}
where $G_{m,\lambda}(x)$ is the cumulative distribution function of
a non-central chi-square variate with $m$ degrees of freedom
and non-centrality parameter $\lambda$. Here,
$\lambda = \phi\tr\{\bm{K}_{22.1}\bm{\epsilon}\bm{\epsilon}^\top\}/2$,
where $\bm{K}_{22.1} = \bm{K}_{\bm{\beta}22} - \bm{K}_{\bm{\beta}21}\bm{K}_{\bm{\beta}11}^{-1}
\bm{K}_{\bm{\beta}12}$ and $\tr(\cdot)$ denotes the trace operator.
The coefficients $b_{ik}$'s ($i=1,2,3,4$ and $k=0,1,2,3$) can be written
in matrix notation, after extensive algebra, as
\begin{align*}
b_{11} &= \frac{\phi}{2}\tr\{(\bm{E} + 2\bm{G})\bm{B}\bm{T}^{(2)}
+ (2\bm{E}-\bm{F}+2\bm{G})\bm{T}^{(3)}+\bm{W}\bm{T}(\bm{C}+2\bm{P})\}\\
&\quad+\frac{1}{2}\tr\{(2\bm{E}-\bm{F}+2\bm{G})\bm{Z}_{1d}\bm{T}+\bm{W}\bm{J}\bm{T}\},
\end{align*}
\[
b_{12} = -\frac{\phi}{6}\tr\{(3\bm{E}-2\bm{F} + 2\bm{G})\bm{T}^{(3)}\},
\quad
b_{13} = 0,
\]
\begin{align*}
b_{21} &= \frac{\phi}{2}\tr\{(\bm{E} + 2\bm{G})\bm{B}\bm{T}^{(2)}
+ (2\bm{E}-\bm{F}+2\bm{G})\bm{T}^{(3)}+\bm{W}\bm{T}(\bm{C}+2\bm{P})\} \\
&\quad + \frac{1}{2}\tr\{(2\bm{E}-\bm{F}+2\bm{G})\bm{Z}_{d}\bm{T}
+ 2(\bm{F}-\bm{E})(\bm{Z}_d - \bm{Z}_{1d})\bm{T}+\bm{W}(\bm{U}\bm{T} + 2\bm{H})\},
\end{align*}
\begin{align*}
b_{22} = \frac{\phi}{2}\tr\{(\bm{F}-\bm{E})\bm{T}^{(3)}+\bm{W}\bm{T}\bm{C}\}
-\frac{1}{2}\tr\{(\bm{F} + 2\bm{G})(\bm{Z}_d - \bm{Z}_{1d})\bm{T}+\bm{W}\bm{T}(\bm{U} - \bm{J}) + 2\bm{W}\bm{H}\},
\end{align*}
\[
b_{23} = -\frac{\phi}{6}\tr\{(\bm{F} + 2\bm{G})\bm{T}^{(3)}+3\bm{W}\bm{T}\bm{C}\},
\]
\begin{align*}
b_{31} &= \frac{\phi}{2}\tr\{(\bm{E} + 2\bm{G})\bm{B}\bm{T}^{(2)}
+ (2\bm{E}-\bm{F}+2\bm{G})\bm{T}^{(3)}+\bm{W}\bm{T}(\bm{C}+2\bm{P})\}\\
&\quad + \frac{1}{2}\tr\{(2\bm{E}-\bm{F}+2\bm{G})\bm{Z}_{1d}\bm{T}
+ (3\bm{E}-2\bm{F} + 2\bm{G})(\bm{Z}_d - \bm{Z}_{1d})\bm{T}+\bm{W}\bm{T}\bm{J}\},
\end{align*}
\[
b_{32} = -\frac{1}{2}\tr\{(3\bm{E}-2\bm{F} + 2\bm{G})(\bm{Z}_d - \bm{Z}_{1d})\bm{T}\},
\quad
b_{33} = -\frac{\phi}{6}\tr\{(3\bm{E}-2\bm{F} + 2\bm{G})\bm{T}^{(3)}\},
\]
\begin{align*}
b_{41} &= \frac{\phi}{2}\tr\{(\bm{E}+2\bm{G})\bm{B}\bm{T}^{(2)}
+ (2\bm{E}-\bm{F}+2\bm{G})\bm{T}^{(3)}+\bm{W}\bm{T}(\bm{C}+2\bm{P})\}\\
&\quad + \frac{1}{4}\tr\{(6\bm{G}-\bm{F} + 4\bm{E})\bm{Z}_{1d}\bm{T}
- (\bm{F} + 2\bm{G})\bm{Z}_{d}\bm{T}+\bm{W}\bm{T}(3\bm{J}-\bm{U}) - 2\bm{W}\bm{H}\},
\end{align*}
\begin{align*}
b_{42} &= -\frac{\phi}{4}\tr\{(2\bm{E}-\bm{F}+2\bm{G})\bm{T}^{(3)}+\bm{W}\bm{T}\bm{C}\}\\
&\quad +\frac{1}{4}\tr\{(\bm{F} + 2\bm{G})(\bm{Z}_d - \bm{Z}_{1d})\bm{T}+\bm{W}\bm{T}(\bm{U}-\bm{J}) + 2\bm{W}\bm{H}\},
\end{align*}
\[
b_{43} = \frac{\phi}{12}\tr\{(\bm{F} + 2\bm{G})\bm{T}^{(3)}+3\bm{W}\bm{T}\bm{C}\},
\]
where
$\bm{U}=\diag\{u_1,\ldots,u_n\}$ with $u_l = \tr\{\bm{X}_{l}^*(\bm{X}^{*\top}\bm{W}\bm{X}^{*})^{-1}\}$,
$\bm{J}=\diag\{j_1,\ldots,j_n\}$ with $j_l = \tr\{\bm{X}_{11l}^*(\bm{X}_{1}^{*\top}\bm{W}\bm{X}_{1}^{*})^{-1}\}$,
$\bm{C}=\diag\{c_1,\ldots,c_n\}$ with $c_{l}=\tr\{\bm{X}_l^*\bm{\epsilon}^*\bm{\epsilon}^{*\top}\}$,
$\bm{P}=\diag\{p_1,\ldots,p_n\}$ with $p_l = \tr\{\bm{X}_{l}^*\bm{\epsilon}^*\bm{\delta}^\top\}$,
$\bm{H}=\diag\{h_1,\ldots,h_n\}$ with $h_l = \phi\tr\{\bm{M}\bm{X}_{l}^*\bm{\epsilon}^*\bm{x}_l^{*\top}\}$,
$\bm{\delta}^\top = (\bm{0}^\top,\bm{\epsilon}^\top)$ and
$\bm{x}_l^{*\top}$ is the $l$th line of $\bm{X}^*$.
The coefficients $b_{i0}$ are obtained from $b_{i0} = -(b_{i1} + b_{i2} + b_{i3})$,
for $i=1,2,3,4$. The $b_{ik}$'s are of order $n^{-1/2}$ and all quantities except $\bm{\epsilon}$
are evaluated under the null hypothesis $\mathcal{H}_0$. The detailed derivation of these expressions
is long and extremely tedious but may be obtained from the authors upon request.

It is interesting to note that the $b_{ik}$'s are functions of the local
derivative matrix and of the (possibly unknown) precision parameter.
These coefficients depend on the second derivative of the (possibly nonlinear)
function $f(\bm{x}_{l};\bm{\beta})$ and
involve the link function and its first and second derivatives.
Unfortunately, they are very difficult to interpret. The
matrices $\bm{C}$, $\bm{H}$, $\bm{J}$, $\bm{P}$ and $\bm{U}$
may be considered the nonlinear contribution of the dispersion model since
they vanish if the regression model is linear.
Obviously, these coefficients depend heavily on the particular dispersion
model under consideration. In particular, these coefficients do not change for the
class of PDMs, since the only difference between PDMs and DMs
is the form of the function $c(\cdot,\cdot)$, which can be
decomposed as $c(y,\phi)=a_1(y)+a_2(\phi)$ for PDMs. 
By replacing $\bm{E}$ by $\bm{F}-\bm{G}$ in these coefficients, we obtain the
nonnull asymptotic distributions of the four statistics in
the class of EFNLMs \citep[see][]{Lemonte2011}. 

Some simplifications in the coefficients $b_{ik}$ ($i=1,2,3,4$ and $k=0,1,2,3$)
can be achieved by examining special cases. For example,
consider the null hypothesis $\mathcal{H}_0:\bm{\beta}=\bm{\beta}_0$ (i.e.~$q=0$) and
an identity link function ($d(\theta_{l})=\theta_{l}$),
which implies that $f_l=-D_{3l}$, $g_l=0$ and $e_l = -D_{2l}'$
($l=1,\ldots,n$). Therefore, the $b_{ik}$'s can be written as
\[
b_{11} = \frac{\phi}{2}\tr\{\bm{E}\bm{B}\bm{T}^{(2)} + (2\bm{E}-\bm{F})\bm{T}^{(3)}+\bm{W}\bm{T}(\bm{C}+2\bm{P})\}
+\frac{1}{2}\tr\{\bm{W}\bm{J}\bm{T}\},
\]
\[
b_{12} = b_{33} = -\frac{\phi}{6}\tr\{(3\bm{E}-2\bm{F})\bm{T}^{(3)}\},
\quad
b_{13} = 0,
\quad
b_{32} = -\frac{1}{2}\tr\{(3\bm{E}-2\bm{F})\bm{Z}_d\bm{T}\},
\]
\begin{align*}
b_{21} &= \frac{\phi}{2}\tr\{\bm{E}\bm{B}\bm{T}^{(2)} + (2\bm{E}-\bm{F})\bm{T}^{(3)}+\bm{W}\bm{T}(\bm{C}+2\bm{P})\} \\
&\quad + \frac{1}{2}\tr\{\bm{F}\bm{Z}_{d}\bm{T}+\bm{W}(\bm{U}\bm{T} + 2\bm{H})\},
\end{align*}
\[
b_{22} = \frac{\phi}{2}\tr\{(\bm{F}-\bm{E})\bm{T}^{(3)}+\bm{W}\bm{T}\bm{C}\}
-\frac{1}{2}\tr\{\bm{F}\bm{Z}_d\bm{T}+\bm{W}\bm{T}(\bm{U} - \bm{J}) + 2\bm{W}\bm{H}\},
\]
\[
b_{23} = -2b_{43} =-\frac{\phi}{6}\tr\{\bm{F}\bm{T}^{(3)}+3\bm{W}\bm{T}\bm{C}\},
\]
\begin{align*}
b_{31} &= \frac{\phi}{2}\tr\{\bm{E}\bm{B}\bm{T}^{(2)} + (2\bm{E}-\bm{F})\bm{T}^{(3)}+\bm{W}\bm{T}(\bm{C}+2\bm{P})\}\\
&\quad + \frac{1}{2}\tr\{(3\bm{E}-2\bm{F})\bm{Z}_d\bm{T}+\bm{W}\bm{T}\bm{J}\},
\end{align*}
\begin{align*}
b_{41} &= \frac{\phi}{2}\tr\{\bm{E}\bm{B}\bm{T}^{(2)} + (2\bm{E}-\bm{F})\bm{T}^{(3)}+\bm{W}\bm{T}(\bm{C}+2\bm{P})\}\\
&\quad + \frac{1}{4}\tr\{-\bm{F}\bm{Z}_{d}\bm{T}+\bm{W}\bm{T}(3\bm{J}-\bm{U}) - 2\bm{W}\bm{H}\},
\end{align*}
\[
b_{42} = -\frac{\phi}{4}\tr\{(2\bm{E}-\bm{F})\bm{T}^{(3)}+\bm{W}\bm{T}\bm{C}\}
+\frac{1}{4}\tr\{\bm{F}\bm{Z}_d\bm{T}+\bm{W}\bm{T}(\bm{U}-\bm{J}) + 2\bm{W}\bm{H}\},
\]
and  $b_{i0} = -(b_{i1} + b_{i2} + b_{i3})$, for $i=1,2,3,4$. For the
log-gamma model, the above coefficients reduce to
\[
b_{11} = \frac{\phi}{2}\tr\{-\bm{F}\bm{T}^{(3)}+\bm{W}\bm{T}(\bm{C}+2\bm{P})\}
+\frac{1}{2}\tr\{\bm{W}\bm{J}\bm{T}\},
\quad
b_{12} = b_{33} = \frac{\phi}{3}\tr\{\bm{F}\bm{T}^{(3)}\},
\quad
b_{13} = 0,
\]
\[
b_{21} = \frac{\phi}{2}\tr\{-\bm{F}\bm{T}^{(3)}+\bm{W}\bm{T}(\bm{C}+2\bm{P})\}
+ \frac{1}{2}\tr\{\bm{F}\bm{Z}_{d}\bm{T}+\bm{W}(\bm{U}\bm{T} + 2\bm{H})\},
\]
\[
b_{22} =\frac{\phi}{2}\tr\{\bm{F}\bm{T}^{(3)}+\bm{W}\bm{T}\bm{C}\}
-\frac{1}{2}\tr\{\bm{F}\bm{Z}_d\bm{T}+\bm{W}\bm{T}(\bm{U} - \bm{J}) + 2\bm{W}\bm{H}\},
\]
\[
b_{23} = -2b_{43} = -\frac{\phi}{6}\tr\{\bm{F}\bm{T}^{(3)}+3\bm{W}\bm{T}\bm{C}\},
\quad
b_{32} = \tr\{\bm{F}\bm{Z}_d\bm{T}\},
\]
\[
b_{31} = \frac{\phi}{2}\tr\{-\bm{F}\bm{T}^{(3)}+\bm{W}\bm{T}(\bm{C}+2\bm{P})\}
+ \frac{1}{2}\tr\{-2\bm{F}\bm{Z}_d\bm{T}+\bm{W}\bm{T}\bm{J}\},
\]
\[
b_{41} = \frac{\phi}{2}\tr\{-\bm{F}\bm{T}^{(3)}+\bm{W}\bm{T}(\bm{C}+2\bm{P})\}
+ \frac{1}{4}\tr\{-\bm{F}\bm{Z}_{d}\bm{T}+\bm{W}\bm{T}(3\bm{J}-\bm{U}) - 2\bm{W}\bm{H}\},
\]
\[
b_{42} = -\frac{\phi}{4}\tr\{-\bm{F}\bm{T}^{(3)}+\bm{W}\bm{T}\bm{C}\}
+\frac{1}{4}\tr\{\bm{F}\bm{Z}_d\bm{T}+\bm{W}\bm{T}(\bm{U}-\bm{J}) + 2\bm{W}\bm{H}\},
\]
Also, for the von Mises model we have
\[
b_{11} = b_{31} = \frac{\phi}{2}\tr\{\bm{W}\bm{T}(\bm{C}+2\bm{P})\}
+\frac{1}{2}\tr\{\bm{W}\bm{J}\bm{T}\},
\quad
b_{12} = b_{13} = b_{32} = b_{33} =0,
\]
\[
b_{21} = \frac{\phi}{2}\tr\{\bm{W}\bm{T}(\bm{C}+2\bm{P})\}
+ \frac{1}{2}\tr\{\bm{W}(\bm{U}\bm{T} + 2\bm{H})\},
\quad
b_{23} = -2b_{43} = -\frac{\phi}{2}\tr\{\bm{W}\bm{T}\bm{C}\},
\]
\[
b_{22}  = -2b_{42} = \frac{\phi}{2}\tr\{\bm{W}\bm{T}\bm{C}\}
-\frac{1}{2}\tr\{\bm{W}\bm{T}(\bm{U} - \bm{J}) + 2\bm{W}\bm{H}\},
\]
\[
b_{41} = \frac{\phi}{2}\tr\{\bm{W}\bm{T}(\bm{C}+2\bm{P})\}
+ \frac{1}{4}\tr\{\bm{W}\bm{T}(3\bm{J}-\bm{U}) - 2\bm{W}\bm{H}\},
\]
Note that for the von Mises linear regression model,
the $b_{ij}$'s above vanish and hence we can write
\[
\Pr(S_{i}\leq x) = G_{p,\lambda}(x) + O(n^{-1}),\qquad i=1,2,3,4.
\]
This is a very interesting
result, which implies that the likelihood ratio, score, Wald
and gradient tests for testing the null hypothesis
$\mathcal{H}_{0}:\bm{\beta}=\bm{\beta}_{0}$
have exactly the same local power up to an error of order $n^{-1}$
when we consider an identity link function.
It should be noticed that  this result also happens for testing the
composite null hypothesis $\mathcal{H}_{0}:\bm{\beta}_{2}=\bm{\beta}_{20}$,
i.e~$\Pr(S_{i}\leq x) = G_{p-q,\lambda}(x) + O(n^{-1})$, for $i=1,2,3,4$.

Now, we present the coefficients that define the nonnull asymptotic distributions
of the likelihood ratio, Wald, score and gradient statistics
for testing the composite null hypothesis
$\mathcal{H}_{0}:\bm{\beta}_{2}=\bm{\beta}_{20}$
in GLMs. We have $t(y_{l},\theta_{l})=y_{l}\theta_{l}-b(\theta_{l})$
and $\mu_{l} = \Es(Y_{l})=\dd b(\theta_{l})/\dd\theta_{l}$.
The class of GLMs is characterized by its variance function
$V_{l} = \dd\mu_{l}/\dd\theta_{l}$, which plays a key role in the study of
its mathematical properties and estimation.
The variance of $Y_{l}$ can be written as var$(Y_{l}) = \phi^{-1}V_{l}$.
For the GLMs, we have $D_{2l}=-V_l^{-1}$ and $D_{3l}=2V_{l}^{-1}(\dd V_l/\dd\mu_l)$
and hence we can rewrite
\[
f_l = \frac{1}{V_l}\frac{\dd \mu_l}{\dd \eta_l}\frac{\dd ^2\mu_l}{\dd \eta_l^2},\quad
g_l = \frac{1}{V_l}\frac{\dd \mu_l}{\dd \eta_l}\frac{\dd ^2\mu_l}{\dd \eta_l^2}
-\frac{1}{V_l^2}\frac{\dd V_l}{\dd \mu_l}\biggl(\frac{\dd \mu_l}{\dd \eta_l}\biggr)^3,
\quad l=1,\ldots,n,
\]
and redefine the matrices $\bm{F}$ and $\bm{G}$ given before.
Additionally, the link function is $d(\mu_{l})=\eta_{l}=\bm{x}_l^\top\bm{\beta}$
with $m=p$. Also, $\bm{\eta}=\bm{X}\bm{\beta}$ with $\bm{X}=(\bm{x}_1,\ldots,\bm{x}_n)^\top$,
i.e.~here $\bm{X}^*=\bm{X}$. Hence, in this class of models we have
\[
b_{11} = \frac{\phi}{2}\tr\{(\bm{F} + \bm{G})\bm{B}\bm{T}^{(2)} + \bm{F}\bm{T}^{(3)}\}
    + \frac{1}{2}\tr\{\bm{F}\bm{Z}_{1d}\bm{T}\},
\quad
b_{12} = b_{33} = -\frac{\phi}{6}\tr\{(\bm{F} - \bm{G})\bm{T}^{(3)}\},
\]
\[
b_{21} = \frac{\phi}{2}\tr\{(\bm{F} + \bm{G})\bm{B}\bm{T}^{(2)} + \bm{F}\bm{T}^{(3)}\}
    + \frac{1}{2}\tr\{\bm{F}\bm{Z}_{d}\bm{T} + 2\bm{G}(\bm{Z}_d - \bm{Z}_{1d})\bm{T}\},
\]
\[
b_{22} = \frac{\phi}{2}\tr\{\bm{G}\bm{T}^{(3)}\}-\frac{1}{2}\tr\{(\bm{F} + 2\bm{G})(\bm{Z}_d - \bm{Z}_{1d})\bm{T}\},
\quad
b_{13} = 0,
\]
\[
b_{23}  = -2b_{43} = -\frac{\phi}{6}\tr\{(\bm{F} + 2\bm{G})\bm{T}^{(3)}\},
\quad
b_{32} = -\frac{1}{2}\tr\{(\bm{F} - \bm{G})(\bm{Z}_d - \bm{Z}_{1d})\bm{T}\},
\]
\[
b_{31} = \frac{\phi}{2}\tr\{(\bm{F} + \bm{G})\bm{B}\bm{T}^{(2)} + \bm{F}\bm{T}^{(3)}\}
    + \frac{1}{2}\tr\{\bm{F}\bm{Z}_{1d}\bm{T} + (\bm{F}-\bm{G})(\bm{Z}_d - \bm{Z}_{1d})\bm{T}\},
\]
\[
b_{41} = \frac{\phi}{2}\tr\{(\bm{F}+\bm{G})\bm{B}\bm{T}^{(2)} + \bm{F}\bm{T}^{(3)}\}
    + \frac{1}{4}\tr\{(3\bm{F} + 2\bm{G})\bm{Z}_{1d}\bm{T} - (\bm{F} + 2\bm{G})\bm{Z}_{d}\bm{T}\},
\]
\[
b_{42} = -\frac{\phi}{4}\tr\{\bm{F}\bm{T}^{(3)}\}
        +\frac{1}{4}\tr\{(\bm{F} + 2\bm{G})(\bm{Z}_d - \bm{Z}_{1d})\bm{T}\},
\]
By considering the identity link function, these coefficients reduce to
\[
b_{11} = \frac{\phi}{2}\tr\{\bm{G}\bm{B}\bm{T}^{(2)}\},
\quad
b_{12} = b_{33} = \frac{\phi}{6}\tr\{\bm{G}\bm{T}^{(3)}\},
\quad
b_{32} = b_{42} = \frac{1}{2}\tr\{\bm{G}(\bm{Z}_d - \bm{Z}_{1d})\bm{T}\},
\]
\[
b_{13} = 0,
\quad
b_{23} = -2b_{43} = -2b_{12}, 
\quad
b_{21} = b_{11} + 2b_{32},
\quad
b_{22} = 3b_{12}- 2b_{32},
\quad
b_{31} = b_{41} = b_{11}-b_{32}.
\]
As expected, the above coefficients vanish for the normal model since
the nonnull distributions of all the four criteria agree with the
$\chi_{p-q,\lambda}^{2}$ distribution.

\section{Power comparisons}\label{power_comparison}

It is known that, to the first order of approximation, the likelihood ratio,
Wald, score and gradient statistics have the same asymptotic
distributional properties either under the null hypothesis or under a sequence of
local alternatives.
On the other hand, up to an error of order $n^{-1}$ the corresponding criteria have the
same size properties but their local powers differ in the $n^{-1/2}$ term. A meaningful
comparison among the criteria can then be performed by comparing
the nonnull asymptotic expansions to order $n^{-1/2}$,
i.e.~ignoring terms or order less than $n^{-1/2}$.

In what follows, we shall compare the local powers of the rival tests
based on the general nonnull asymptotic expansions derived in Section \ref{main_result}
for testing the null hypothesis $\mathcal{H}_0:\bm{\beta}_{2}=\bm{\beta}_{20}$
in the class of DMs. Let $\Pi_{i}$ be the power function, up to order $n^{-1/2}$, of the test that uses
the statistic $S_{i}$, for $i=1,2,3,4$. We have
\begin{equation}\label{diff_power}
\Pi_{i} - \Pi_{j} = \sum_{k=0}^{3}(b_{jk} - b_{ik})G_{p-q+2k,\lambda}(x),
\end{equation}
for $i\neq j$. It is well known that
\begin{equation}\label{diff_G}
G_{m,\lambda}(x) - G_{m+2,\lambda}(x) = 2g_{m+2,\lambda}(x),
\end{equation}
where $g_{\nu,\lambda}(x)$ is the probability density
function of a non-central chi-square random variable
with $\nu$ degrees of freedom and non-centrality parameter $\lambda$.
From~(\ref{diff_power}) and (\ref{diff_G}) we have after some algebra
\begin{align}\label{eq:power}
\begin{split}
\Pi_{1} - \Pi_{4} &= k_{1}g_{p-q+4,\lambda}(x) + k_{2}g_{p-q+6,\lambda}(x),\quad
\Pi_{2} - \Pi_{4}  = k_{3}g_{p-q+4,\lambda}(x) + k_{4}g_{p-q+6,\lambda}(x),\\
\Pi_{3} - \Pi_{4} &= k_{5}g_{p-q+4,\lambda}(x) + k_{6}g_{p-q+6,\lambda}(x),\quad
\Pi_{1} - \Pi_{2}  = k_{7}g_{p-q+4,\lambda}(x) + k_{8}g_{p-q+6,\lambda}(x),\\
\Pi_{1} - \Pi_{3} &= k_{9}g_{p-q+4,\lambda}(x) + k_{10}g_{p-q+6,\lambda}(x),\quad
\Pi_{2} - \Pi_{3}  = k_{11}g_{p-q+4,\lambda}(x) + k_{12}g_{p-q+6,\lambda}(x),
\end{split}
\end{align}
where
\[
k_1 = -\frac{1}{2}\tr\{(\bm{F}+2\bm{G})(\bm{Z}_d - \bm{Z}_{1d})\bm{T}\}
+ \frac{1}{2}\tr\{\bm{W}\bm{T}(\bm{J}-\bm{U})-2\bm{W}\bm{H}\},
\]
\[
k_2 = -\frac{\phi}{6}\tr\{(\bm{F}+2\bm{G})\bm{T}^{(3)}\} - \frac{\phi}{2}\tr\{\bm{W}\bm{T}\bm{C}\},
\quad
k_3 = 3k_1,
\quad
k_4 = 3k_2,
\]
\[
k_5 = k_1 - \tr\{(3\bm{E}-2\bm{F}+2\bm{G})(\bm{Z}_d - \bm{Z}_{1d})\bm{T}\},
\]
\[
k_6 = - \frac{\phi}{2}\tr\{(2\bm{E}-\bm{F}+2\bm{G})\bm{T}^{(3)}\}- \frac{\phi}{2}\tr\{\bm{W}\bm{T}\bm{C}\},
\]
\[
k_7 = -2k_1,
\quad
k_8 = -2k_2,
\quad
k_9 = k_1 - k_5,
\quad
k_{10} = \frac{\phi}{3}\tr\{(3\bm{E}-2\bm{F}+2\bm{G})\bm{T}^{(3)}\},
\]
\[
k_{11} = -3\tr\{(\bm{F}-\bm{E})(\bm{Z}_d- \bm{Z}_{1d})\bm{T}\}
-\tr\{\bm{W}\bm{T}(\bm{U}-\bm{J})+2\bm{W}\bm{H}\},
\]
\[
k_{12} = -\phi\tr\{(\bm{F}-\bm{E})\bm{T}^{(3)}\}
-\phi\tr\{\bm{W}\bm{T}\bm{C}\}.
\]

For proper dispersion models, the above expressions are the same.
Replacing $\bm{E}$ by $\bm{F}-\bm{G}$
we obtain these quantities for exponential family nonlinear models.
From equations~(\ref{eq:power}) we have $\Pi_{1}>\Pi_{3}$ if $k_{9}\geq 0$ and $k_{10}\geq 0$
with $k_{9}+k_{10}>0$, and if $k_{9}\leq 0$ and $k_{10}\leq 0$ with $k_{9}+k_{10}<0$, we have
$\Pi_{1}<\Pi_{3}$. Also, $\Pi_{1}=\Pi_{3}$ if $k_{9}=k_{10}=0$, i.e.~$\bm{F}=\bm{G}$
and $\bm{E}=\bm{0}$, which occurs only for von Mises and normal models with any link function.
Additionally, equations (\ref{eq:power})
show that with the exception of the likelihood ratio
and score tests, is not possible to have any other
equality among the power functions in the class of DMs for
testing the null hypothesis $\mathcal{H}_0:\bm{\beta}_{2}=\bm{\beta}_{20}$.
The reason is that $\bm{C}$, $\bm{H}$, $\bm{J}$ and $\bm{U}$,
which may be considered as the nonlinear contribution of the dispersion model,
vanish only for linear regression models.
It implies that only strict inequality holds for any other power comparison
among the power functions of the tests that are based on the
statistics $S_1$, $S_2$, $S_3$ and $S_4$. For example,
from (\ref{eq:power}) we have $\Pi_{1}>\Pi_{4}$ ($\Pi_{1}<\Pi_{4}$)
if $k_{1}\geq 0$ and $k_{2}\geq 0$ with $k_{1}+k_{2}>0$
(if $k_{1}\leq 0$ and $k_{2}\leq 0$ with $k_{1}+k_{2}<0$), and so on.

We now move to the class of GLMs, in which 
$\bm{C}=\bm{H}=\bm{J}=\bm{P}=\bm{U}=\bm{0}$. By using the coefficients derived
for this class of models in Section \ref{main_result}, the quantities that define
equation (\ref{eq:power}) reduce to
\[
k_1 = -\frac{1}{2}\tr\{(\bm{F}+2\bm{G})(\bm{Z}_d - \bm{Z}_{1d})\bm{T}\},
\quad
k_2 = -\frac{\phi}{6}\tr\{(\bm{F}+2\bm{G})\bm{T}^{(3)}\},
\quad
k_3 = 3k_1,
\]
\[
k_5 = k_1 - \tr\{(\bm{F}-\bm{G})(\bm{Z}_d - \bm{Z}_{1d})\bm{T}\},
\quad
k_6 = - \frac{\phi}{2}\tr\{\bm{F}\bm{T}^{(3)}\},
\quad
k_4 = 3k_2,
\]
\[
k_7 = -2k_1,
\quad
k_8 = -2k_2,
\quad
k_9 = k_1 - k_5,
\quad
k_{10} = \frac{\phi}{3}\tr\{(\bm{F}-\bm{G})\bm{T}^{(3)}\},
\]
\[
k_{11} = -3\tr\{\bm{G}(\bm{Z}_d- \bm{Z}_{1d})\bm{T}\},
\quad
k_{12} = -\phi\tr\{\bm{G}\bm{T}^{(3)}\}.
\]
For GLMs with canonical link ($\bm{G}=\bm{0}$), we have
$k_{11}=k_{12}=0$ and hence $\Pi_{2}=\Pi_{3}$.
It is possible to show that $\Pi_{1}=\Pi_{2}=\Pi_{4}$ if $\bm{F}=-2\bm{G}$, that is
\[
\frac{\dd^2\mu_{l}}{\dd\eta_l^2} = \frac{2}{3V_l}\biggl(\frac{\dd\mu_{l}}{\dd\eta_l}\biggr)^2,
\quad l=1,\ldots,n.
\]
The GLMs for which this equality holds have the link function
defined by $\eta_{l}=\int V_l^{-3/2}\dd\mu_{l}$ ($l=1,\ldots,n$). For the gamma
model this function is $\eta_{l}=\mu_{l}^{-1/3}$ ($l=1,\ldots,n$).
Additionally, we have that $\Pi_{3}=\Pi_{4}$ for any GLM with identity link function,
i.e.~$\bm{F}=\bm{0}$. Also, $\Pi_{1}=\Pi_{3}$ if $k_{9}=k_{10}=0$,
i.e.~$\bm{F}=\bm{G}$, which occurs only for normal models
with any link. Finally, the equality $\Pi_{1}=\Pi_{2}=\Pi_{3}=\Pi_{4}$ holds
only for normal models with identity link function.

We can conclude that there is no uniform superiority
of one test with respect to the others for
testing the null hypothesis $\mathcal{H}_{0}:\bm{\beta}_{2}=\bm{\beta}_{20}$
in the class of DMs. Hence, if the sample size is large, all tests
could be recommended, since their type I error probabilities
do not significantly deviate from the true nominal level
and their local powers are approximately equal.
The natural question is how these tests perform when the sample size is
small or of moderate size, and which one is the most reliable.
In Section \ref{simulations}, we shall use Monte Carlo simulations
to shed some light on this issue.

\section{Tests for the precision parameter}\label{test_phi}

In this section we derive asymptotic expansions for the nonnull
distribution of the four statistics for testing the precision
parameter $\phi$ in DMs. We are interested in testing the
null hypothesis $\mathcal{H}_{0}:\phi=\phi_{0}$ against
a two-sided alternative hypothesis
$\mathcal{H}_{1}:\phi\neq\phi_{0}$, where $\phi_{0}$ is a positive specified value
for $\phi$. Here, $\bm{\beta}$ acts as a nuisance parameter. The likelihood
ratio, Wald, score and gradient
statistics are expressed as follows:
\[
S_1 = \sum_{l=1}^{n}\{(\widehat{\phi}-\phi_{0})t(y_l,\widehat{\theta}_{l})
+ c(y_l,\widehat{\phi}) - c(y_l,\phi_{0})\},
\quad
S_{2} = (\widehat{\phi} - \phi_{0})^2\{-\alpha_{2}(\widehat{\phi})\},
\]
\[
S_{3} = \{-\alpha_{2}(\phi_{0})\}^{-1}\Biggl[\sum_{l=1}^n\{t(y_l,\widehat{\theta}_{l})+c^{(1)}(y_l,\phi_{0})\}\Biggr]^2,
\quad
S_{4} = (\widehat{\phi} - \phi_{0})\sum_{l=1}^n\{t(y_l,\widehat{\theta}_{l})+c^{(1)}(y_l,\phi_{0})\}.
\]
For PDMs, these statistics can be expressed as
\[
S_{1} = 2n\{a_{2}(\widehat{\phi}) - a_{2}(\phi_{0})
- (\widehat{\phi} - \phi_{0})a_2^{(1)}(\widehat{\phi})\},
\quad
S_{2} = -n(\widehat{\phi} - \phi_{0})^2a_{2}^{(2)}(\widehat{\phi}),
\]
\[
S_{3} = -\frac{n\{a_{2}^{(1)}(\widehat{\phi}) - a_{2}^{(1)}(\phi_{0})\}^2}{a_{2}^{(2)}(\phi_{0})},
\quad
S_{4} = n\{a_{2}^{(1)}(\phi_{0}) - a_{2}^{(1)}(\widehat{\phi})\}(\widehat{\phi} - \phi_{0}).
\]
For example, for the von Mises model $a_2(\phi) = -\log\{I_0(\phi)\}$. Also,
$a_2^{(1)}(\phi) = -r(\phi)$ and $a_2^{(2)}(\phi) = r(\phi)^2 + r(\phi)/\phi-1$, where
$r(\phi) = I_1(\phi)/I_0(\phi)$. Thus, we can write
\[
S_{1} = 2n[\log\{I_0(\phi_0)/I_0(\widehat{\phi})\}
+ (\widehat{\phi} - \phi_{0})r(\widehat{\phi})],
\quad
S_{2} = -n(\widehat{\phi} - \phi_{0})^2\{r(\widehat{\phi})^2 + r(\widehat{\phi})/\widehat{\phi}-1\},
\]
\[
S_{3} = -\frac{n\{r(\phi_{0})-r(\widehat{\phi})\}^2}{r(\phi_{0})^2 + r(\phi_{0})/\widehat{\phi_{0}}-1},
\quad
S_{4} = n\{r(\widehat{\phi})-r(\phi_{0})\}(\widehat{\phi} - \phi_{0}).
\]
Also, for normal and inverse Gaussian models we have $a_{2}(\phi) = \log(\phi)/2$. Hence
\[
S_{1} = 2n\biggl\{\log\biggl(\frac{\widehat{\phi}}{\phi_{0}}\biggr)
- \biggl(\frac{\widehat{\phi} - \phi_{0}}{\widehat{\phi}}\biggr)\biggr\},
\quad
S_{2} = S_{3} = \frac{n}{2}\biggl\{\frac{\widehat{\phi} - \phi_{0}}{\widehat{\phi}}\biggr\}^2,
\quad
S_{4} = \frac{n}{2}\biggl\{\frac{\widehat{\phi} - \phi_{0}}{\phi_{0}}
- \frac{\widehat{\phi} - \phi_{0}}{\widehat{\phi}}\biggr\}.
\]
We have $a_2(\phi)=\phi\log(\phi) - \log\{\Gamma(\phi)\}$ for the gamma model
and therefore these statistics reduce to
\[
S_1 = 2n\biggl\{\phi_{0}\log\biggl(\frac{\widehat{\phi}}{\phi_{0}}\biggr)
-\log\biggl(\frac{\Gamma(\widehat{\phi})}{\Gamma(\phi_{0})}\biggr)
-(\widehat{\phi} - \phi_{0})(1-\psi(\widehat{\phi}))\biggr\},
\]
\[
S_2 = n\{\widehat{\phi}\psi'(\widehat{\phi}) - 1\}
\frac{(\widehat{\phi} - \phi_{0})^2}{\widehat{\phi}},
\quad
S_3 = \frac{n\phi_{0}\{\log(\widehat{\phi}/\phi_{0}) - (\psi(\widehat{\phi}) - \psi(\phi_{0}))\}}
{\phi_{0}\psi'(\phi_{0}) - 1}
\]
and
\[
S_{4} = n(\widehat{\phi}-\phi_{0})\biggl\{\log\biggl(\frac{\widehat{\phi}}{\phi_{0}}\biggr)
+\psi(\widehat{\phi}) - \psi(\phi_{0})\biggr\},
\]
where $\Gamma(\cdot)$, $\psi(\cdot)$ and $\psi'(\cdot)$ are the gamma, digamma and
trigamma functions, respectively.

The nonnull asymptotic distributions of $S_1$, $S_2$, $S_3$ and
$S_4$ for testing $\mathcal{H}_{0}:\phi=\phi_{0}$ in DMs under the
local alternative $\mathcal{H}_{1n}:\phi=\phi_{0}+\epsilon$,
where $\epsilon=\phi-\phi_{0}$ is assumed to be $O(n^{-1/2})$, is
\[
\Pr(S_{i}\leq x) = G_{1,\lambda}(x) + \sum_{k=0}^{3}b_{ik}G_{1+2k,\lambda}(x) + O(n^{-1}),
\qquad i=1,2,3,4.
\]
The noncentrality parameter is given by $\lambda = -\alpha_{2}\epsilon^2$
and the the coefficients $b_{ik}$'s can be written as
\[
b_{11} = \frac{(\alpha_2'-\alpha_3)\epsilon^3}{2}+\frac{p\epsilon}{2\phi},
\quad
b_{12} =  \frac{(2\alpha_3-3\alpha_2')\epsilon^3}{6},
\quad
b_{13} = 0,
\]
\[
b_{21} = \frac{(\alpha_2'-\alpha_3)\epsilon^3}{2}-\frac{\alpha_{3}\epsilon}{2\alpha_2}+\frac{p\epsilon}{2\phi},
\quad
b_{22} = -\frac{(\alpha_2'-\alpha_3)\epsilon^3}{2}+\frac{\alpha_{3}\epsilon}{2\alpha_2},
\quad
b_{23} = -\frac{\alpha_3\epsilon^3}{6},
\]
\[
b_{31} = \frac{(\alpha_2'-\alpha_3)\epsilon^3}{2} + \frac{(2\alpha_3-3\alpha_2')\epsilon}{2\alpha_2}
+\frac{p\epsilon}{2\phi},
\quad
b_{32} = -\frac{(2\alpha_3-3\alpha_2')\epsilon}{2\alpha_2},
\quad
b_{33} = \frac{(2\alpha_3-3\alpha_2')\epsilon^3}{6},
\]
\[
b_{41} = \frac{(\alpha_2'-\alpha_3)\epsilon^3}{2} + \frac{\alpha_3\epsilon}{4\alpha_2}
+\frac{p\epsilon}{2\phi},
\quad
b_{42} = -\frac{(2\alpha_2'-\alpha_3)\epsilon^3}{4} - \frac{\alpha_3\epsilon}{4\alpha_2},
\quad
b_{43} = \frac{\alpha_3\epsilon^3}{12},
\]
with $b_{i0} = -(b_{i1}+b_{i2}+b_{i3})$, for $i=1,2,3,4$.
It should be noticed that the above expressions depend on the parameter
$\phi$ and depend on the local derivative matrix $\bm{X}^*$ only through its rank $p$.
Since $\alpha_2'=\alpha_3=na_{2}^{(3)}(\phi)$ for PDMs, these coefficients reduce to
\[
b_{11} =  \frac{p\epsilon}{2\phi},
\quad
b_{12} = b_{23} = b_{33} = -\frac{na_{2}^{(3)}(\phi)\epsilon^3}{6},
\quad
b_{13} = 0,
\quad
b_{21} = b_{31} = \frac{p\epsilon}{2\phi} - \frac{a_{2}^{(3)}(\phi)\epsilon}{2a_{2}^{(2)}(\phi)},
\]
\[
b_{22} = b_{32} = b_{11}-b_{21}, 
\quad
b_{41} = b_{11} + \frac{1}{2}(b_{11}-b_{21}),
\quad
b_{42} = -\frac{1}{2}(b_{11}-b_{21}-3b_{12}),
\quad
b_{43} = -\frac{b_{12}}{2},
\]
with $b_{i0} = -(b_{i1}+b_{i2}+b_{i3})$, for $i=1,2,3,4$.
These coefficients do not change for the class of GLMs.

In what follows, we present an analytical comparison
among the local powers of the four tests for testing
the null hypothesis $\mathcal{H}_{0}:\phi=\phi_{0}$.
We have
\[
\Pi_{i} - \Pi_{j} = \sum_{k=0}^{3}(b_{jk} - b_{ik})G_{1+2k,\lambda}(x).
\]
After some algebra, we can write
\[
\Pi_{1}-\Pi_{2} = -\frac{\alpha_3\epsilon}{\alpha_2} g_{5,\lambda}(x)
+ \frac{\alpha_3\epsilon^3}{3} g_{7,\lambda}(x),
\]
\[
\Pi_{1}-\Pi_{3} = \frac{(2\alpha_3-3\alpha_2')\epsilon}{\alpha_2} g_{5,\lambda}(x)
-\frac{(2\alpha_3-3\alpha_2')\epsilon^3}{3} g_{7,\lambda}(x),
\]
\[
\Pi_{1}-\Pi_{4} = \frac{\alpha_3\epsilon}{2\alpha_2} g_{5,\lambda}(x)
- \frac{\alpha_3\epsilon^3}{6} g_{7,\lambda}(x),
\]
\[
\Pi_{2}-\Pi_{3} = \frac{3(\alpha_3-\alpha_2')\epsilon}{\alpha_2} g_{5,\lambda}(x)
-(\alpha_3-\alpha_2')\epsilon^3 g_{7,\lambda}(x),
\]
\[
\Pi_{2}-\Pi_{4} = \frac{3\alpha_3\epsilon}{2\alpha_2} g_{5,\lambda}(x)
- \frac{\alpha_3\epsilon^3}{2} g_{7,\lambda}(x),
\]
\[
\Pi_{3}-\Pi_{4} = -\frac{3(\alpha_3-2\alpha_2')\epsilon}{\alpha_2} g_{5,\lambda}(x)
+\frac{(\alpha_3-2\alpha_2')}{2}\epsilon^3 g_{7,\lambda}(x).
\]
From the above expressions, we can obtain the following general conclusions.
By assuming $\phi>\phi_{0}$ (opposite inequalities hold if $\phi<\phi_{0}$),
we have that $\Pi_{3} < \Pi_{2} < \Pi_{1} < \Pi_{4}$ if $\alpha_3>0$ with
$\alpha_2'>\alpha_3$. Also, $\Pi_{2} = \Pi_{3} < \Pi_{1} < \Pi_{4}$ if
$\alpha_2'=\alpha_3>0$. For example, for normal and inverse Gaussian models
we have $a_2(\phi) = \log(\phi)/2$,
which implies that $a_2^{(1)}(\phi) = 1/(2\phi)$,
$a_2^{(2)}(\phi) = -1/(2\phi^2)$ and $a_2^{(3)}(\phi) = 1/\phi^3$.
Since $\alpha_2'=\alpha_3=n/\phi^3>0$,
we arrive at the following inequalities:
$\Pi_{2} = \Pi_{3} < \Pi_{1} < \Pi_{4}$  if $\phi>\phi_{0}$,  and
$\Pi_{2} = \Pi_{3} > \Pi_{1} > \Pi_{4}$ if $\phi<\phi_{0}$.

\section{Monte Carlo simulation}\label{simulations}

In this section we conduct Monte Carlo simulations in order
to compare the performance of the likelihood ratio, Wald,
score and gradient tests in small- and moderate-sized samples.

We consider the von Mises regression model, which is quite useful
for modeling circular data; see \cite{Fisher1993} and \cite{MardiaJupp2000}.
Here,
\[
\pi(y;\theta,\phi)=\frac{\exp\{\phi\cos(y-\theta)\}}{2\pi I_0(\phi)},
\quad y\in(-\pi,\pi),
\]
where $\theta\in(-\pi,\pi)$ and $\phi>0$. This density function
is symmetric around $y=\theta$, which is the mode and the circular mean of
the distribution. Also, $\phi$ is a precision parameter in the sense that
the larger the value of $\phi$ the more concentrated the density
function around $\theta$. It is evident the density function above
is a special case of (\ref{dens1}) with $t(y,\theta)=\cos(y-\theta)$ and
$c(y,\phi)=-\log(I_0(\phi))$.

We assume that
\[
\tan(\theta_{l}/2) = \eta_{l} = \beta_{1}x_{i1} + \beta_{2}x_{i2} + \cdots + \beta_{p}x_{ip},
\]
where $x_{i1} = 1$ and $\theta_{l}=2\arctan(\eta_{l})$, $l = 1, \ldots, n$.
The covariate values were selected as random draws from the $\mathcal{U}(0,1)$
distribution and for fixed $n$ those values were kept constant
throughout the experiment. The number of Monte Carlo replications was 10,000, the nominal levels
of the tests were $\gamma$ = 10\%, 5\% and 1\%, and all simulations were carried
out using the {\tt Ox} matrix programming language \citep{DcK2007}.
{\tt Ox} is freely distributed for academic
purposes and available at http://www.doornik.com.

First, the null hypothesis is $\mathcal{H}_{0}:\beta_{p-1} =\beta_{p} = 0$,
which is tested against a two-sided alternative. The sample size is $n=50$,
$\phi = 1.5, 2.5, 4$ and $p = 3, 4, \ldots, 8$. The values of the
response were generated using $\beta_{1} = \cdots=\beta_{p-2} = 1$.
The null rejection rates of the four tests are presented in Table~\ref{tab1}.
It is clear that the likelihood ratio ($S_1$) and Wald ($S_2$) tests are
markedly liberal, more so as the number of regressors increases.
The score ($S_3$) and gradient ($S_4$) tests  are also liberal in most of
the cases, but much less size distorted than the likelihood ratio
and Wald tests in all cases. For instance, when $\phi=2.5$,
$p=4$ and $\gamma = 5\%$, the rejection rates
are 7.05\% ($S_1$), 8.28\% ($S_2$),  5.15\% ($S_3$) and 6.30\% ($S_4$).
We note that the score test is much less liberal than the
likelihood ratio and Wald tests and slightly less liberal than the
gradient test. Additionally, the Wald test is much more liberal than
the other tests. Note that as $\phi$ increases the tests become less size
distorted, as expected, since the
von Mises distribution approaches a normal distribution as $\phi$
increases.
\begin{table}[!htp]
{\footnotesize
\begin{center}
\caption{Null rejection rates (\%); $\phi$ = 1.5, 2.5 and 4, with $n = 50$.}\label{tab2}
\begin{tabular}{c c c c c| c c c c| c c c c  }\hline
  & \multicolumn{12}{c}{$\phi = 1.5$} \\\cline{2-13}
  & \multicolumn{4}{c|}{$\gamma = 10\%$} & \multicolumn{4}{c|}{$\gamma = 5\%$}
  & \multicolumn{4}{c}{$\gamma = 1\%$}\\\cline{2-13}
  $p$   & $S_1$  & $S_2$  & $S_3$ & $S_4$
        & $S_1$  & $S_2$  & $S_3$ & $S_4$
        & $S_1$  & $S_2$  & $S_3$ & $S_4$\\\hline
    3   & 13.31 & 15.42 & 10.12 & 10.42 &  6.90 &  9.93 & 4.65 & 5.04 & 1.75 & 4.13 & 0.79 & 1.20 \\
    4   & 14.48 & 16.31 & 10.26 & 12.49 &  7.75 & 10.86 & 4.83 & 6.83 & 1.93 & 4.62 & 0.59 & 2.08 \\
    5   & 16.65 & 19.34 & 10.92 & 12.46 &  9.55 & 12.36 & 5.05 & 6.62 & 2.67 & 4.87 & 0.84 & 1.83 \\
    6   & 19.04 & 21.93 & 11.94 & 14.81 & 11.78 & 15.00 & 5.90 & 8.26 & 3.62 & 6.50 & 1.03 & 2.40 \\
    7   & 22.09 & 26.39 & 12.44 & 15.94 & 13.71 & 18.12 & 6.12 & 8.87 & 4.27 & 7.67 & 1.27 & 2.21 \\
    8   & 24.16 & 26.58 & 13.03 & 17.66 & 15.87 & 17.42 & 6.63 & 9.82 & 5.23 & 6.82 & 1.39 & 2.76 \\ \hline

  & \multicolumn{12}{c}{$\phi = 2.5$} \\\cline{2-13}
  & \multicolumn{4}{c|}{$\gamma = 10\%$} & \multicolumn{4}{c|}{$\gamma = 5\%$}
  & \multicolumn{4}{c}{$\gamma = 1\%$}\\\cline{2-13}
  $p$   & $S_1$  & $S_2$  & $S_3$ & $S_4$
        & $S_1$  & $S_2$  & $S_3$ & $S_4$
        & $S_1$  & $S_2$  & $S_3$ & $S_4$\\\hline
    3   & 12.02 & 12.96 & 10.56 & 10.50 &  6.21 &  7.35 & 5.17 & 5.29 & 1.39 & 2.31 & 0.78 & 1.04 \\
    4   & 12.97 & 13.66 & 11.05 & 11.77 &  7.05 &  8.28 & 5.15 & 6.30 & 1.73 & 3.05 & 0.90 & 1.52 \\
    5   & 14.28 & 16.38 & 10.97 & 11.68 &  7.96 & 10.31 & 4.94 & 6.25 & 2.11 & 4.28 & 0.85 & 1.65 \\
    6   & 14.83 & 15.33 & 11.90 & 13.02 &  8.36 &  9.82 & 5.71 & 7.27 & 2.09 & 3.85 & 1.01 & 1.80 \\
    7   & 15.93 & 18.00 & 12.60 & 13.87 &  9.20 & 11.30 & 6.66 & 7.60 & 2.72 & 3.71 & 1.53 & 1.87 \\
    8   & 18.12 & 19.53 & 13.45 & 16.12 & 11.16 & 12.29 & 7.02 & 9.38 & 3.31 & 4.79 & 1.55 & 2.68 \\ \hline

  & \multicolumn{12}{c}{$\phi = 4$} \\\cline{2-13}
  & \multicolumn{4}{c|}{$\gamma = 10\%$} & \multicolumn{4}{c|}{$\gamma = 5\%$}
  & \multicolumn{4}{c}{$\gamma = 1\%$}\\\cline{2-13}
  $p$   & $S_1$  & $S_2$  & $S_3$ & $S_4$
        & $S_1$  & $S_2$  & $S_3$ & $S_4$
        & $S_1$  & $S_2$  & $S_3$ & $S_4$\\\hline
    3   & 11.99 & 12.59 & 10.72 & 10.81 & 6.32 &  7.19 & 5.02 & 5.25 & 1.37 & 2.20 & 0.82 & 1.12 \\
    4   & 13.15 & 14.48 & 11.49 & 11.74 & 7.19 &  8.66 & 5.50 & 5.83 & 1.67 & 2.89 & 0.84 & 1.13 \\
    5   & 13.59 & 13.67 & 11.87 & 12.26 & 7.21 &  7.64 & 5.72 & 6.25 & 1.68 & 2.50 & 0.96 & 1.35 \\
    6   & 14.08 & 15.60 & 11.85 & 12.65 & 7.57 &  9.04 & 5.88 & 6.30 & 1.73 & 2.88 & 1.00 & 1.21 \\
    7   & 15.16 & 16.42 & 12.79 & 13.52 & 8.34 &  9.55 & 6.42 & 7.03 & 2.28 & 3.16 & 1.43 & 1.71 \\
    8   & 16.14 & 17.36 & 13.53 & 14.57 & 9.28 & 10.31 & 7.13 & 7.84 & 2.42 & 2.96 & 1.28 & 1.61 \\ \hline
\end{tabular}
\end{center}    }
\end{table}

Table~\ref{tab3} reports results for $\phi = 3$, $p=4$ and sample sizes ranging
from 20 to 150. As expected, the null rejection rates of all the tests approach the
corresponding nominal levels as the sample size grows. Again,
the score and gradient tests present the best performances.
In Table \ref{tab4} we present the first two moments
of $S_1$, $S_2$, $S_3$ and $S_4$ and the corresponding moments
of the limiting $\chi^2$ distribution. Note that the
gradient and score statistics present a good
agreement between the true moments (obtained by simulation)
and the moments of the limiting distribution.
\begin{table}[!htp]
{\footnotesize
\begin{center}
\caption{Null rejection rates (\%); $\phi = 3$,
        $p = 4$ and different sample sizes.}\label{tab3}
\begin{tabular}{c c c c c| c c c c| c c c c  }\hline
  & \multicolumn{4}{c|}{$\gamma = 10\%$} & \multicolumn{4}{c|}{$\gamma = 5\%$}
  & \multicolumn{4}{c}{$\gamma = 1\%$}\\\cline{2-13}
  $n$   & $S_1$  & $S_2$  & $S_3$ & $S_4$
        & $S_1$  & $S_2$  & $S_3$ & $S_4$
        & $S_1$  & $S_2$  & $S_3$ & $S_4$\\\hline
   20   & 17.33 & 19.18 & 13.71 & 13.89 & 10.50 & 11.95 & 6.92 & 7.04 & 3.33 & 4.38 & 1.16 & 1.14 \\
   30   & 15.04 & 16.33 & 11.65 & 12.76 &  8.29 & 10.19 & 5.10 & 6.66 & 2.05 & 4.14 & 0.75 & 1.50 \\
   40   & 13.49 & 15.23 & 11.44 & 11.44 &  7.56 &  9.43 & 5.72 & 5.96 & 1.81 & 3.07 & 0.92 & 1.18 \\
   50   & 12.51 & 13.78 & 10.77 & 11.05 &  6.65 &  7.79 & 5.40 & 5.59 & 1.66 & 2.31 & 1.02 & 1.25 \\
   70   & 12.01 & 12.46 & 11.00 & 11.17 &  6.20 &  6.90 & 5.41 & 5.58 & 1.48 & 2.18 & 1.12 & 1.28 \\
  100   & 11.30 & 12.13 & 10.74 & 10.69 &  5.86 &  6.65 & 4.92 & 5.44 & 1.22 & 2.04 & 0.94 & 1.07 \\
  150   & 10.51 & 11.01 & 10.02 & 10.10 &  5.05 &  6.03 & 4.59 & 4.63 & 1.08 & 1.66 & 0.94 & 0.95 \\ \hline
\end{tabular}
\end{center}  }
\end{table}
\begin{table}[!htp]
{\footnotesize
\begin{center}
\caption{Moments; $\phi = 2$, $n = 35$, $p=4$.}\label{tab4}
\begin{tabular}{l c c c c c}\hline
  & $S_1$ & $S_2$ & $S_3$ & $S_4$ & $\chi_{2}^{2}$ \\\hline
Mean     & 2.50  & 2.68 & 2.16 & 2.23  &  2.0 \\
Variance & 6.23  & 8.73 & 4.14 & 4.63  &  4.0 \\ \hline
\end{tabular}
\end{center}   }
\end{table}

We also performed Monte Carlo simulations considering hypothesis testing on $\phi$.
To save space, the results are not shown. The score and gradient tests exhibited superior
behaviour than the likelihood ratio and Wald tests. For example, when $n=35$, $p=3$,
$\gamma = 10\%$ and $\mathcal{H}_{0}:\phi = 2$, we
obtained the following null rejection rates: 13.23\%
($S_1$), 14.75\% ($S_2$), 10.61\% ($S_3$) and 9.97\% ($S_4$).
Again, the best performing tests are the score and gradient tests.

Overall, in small to moderate-sized samples the best performing tests are the
score and the gradient tests. They are less size distorted than the other two.
Hence, these tests may be recommended for testing hypotheses on the regression parameters
in the von Mises regression model.
The gradient test has a slight advantage over the score test because
the gradient statistic is simpler to calculate than the
score statistic for testing a subset of regression parameters.
In particular, no matrix needs to be inverted; see Section \ref{DMs}.

\section{Application}\label{applications}

In this section we shall illustrate an application of 
the likelihood ratio, Wald, score and gradient tests in a real data set.
We consider the data described in \cite{FL1992} regarding the distance traveled by
31 small blue periwinkles ({\it Nodilittorina unifasciata}) after they have moved down-shore
from the height at which they normally live. Following \cite{FL1992} we assume
a von Mises distribution for the animals' path, but with the assumption of
constant dispersion and link function
\[
\tan(\theta_{l}/2) = \beta_{1} + \beta_{2}x_l,\quad l=1,\ldots,31,
\]
where $\theta_{l}=2\arctan(\beta_{1} + \beta_{2}x_l)$ denotes the mean direction for a given distance
moved $x_{l}$ (cm). These data have been previously analysed by
\cite{Paula1996} and \cite{SouzaPaula2002} with emphasis on local influence and
residual analysis, respectively. The angular responses were transformed
to the range $(-\pi,\pi)$. The maximum likelihood estimates of the parameters
(asymptotic standard errors in parentheses) are: $\widehat{\beta}_{1} = -0.323\,(0.151)$,
$\widehat{\beta}_{2} = -0.013\,(0.004)$ and $\widehat{\phi}=3.265\,(0.726)$.
The values of the likelihood ratio ($S_1$), Wald ($S_2$), score ($S_3$) and
gradient ($S_4$) statistics for testing the null hypothesis
$\mathcal{H}_{0}:\beta_{2}=0$ are 9.526\,($p$-value: 0.002),
11.031\,($p$-value: 0.001), 7.126\,($p$-value: 0.008) and 8.280\,($p$-value: 0.004),
respectively. At any usual significance level, all tests lead to the same
conclusion, i.e.~the null hypothesis should be rejected.

Now, we consider different values for $\beta_{20}$ and we wish to test
$\mathcal{H}_0:\beta_{2}=\beta_{20}$ against $\mathcal{H}_1:\beta_{2}\neq\beta_{20}$.
Table \ref{tab_new} lists the observed values of the different test statistics and the corresponding
$p$-values for $\beta_{20} = -0.026, -0.024, -0.022, -0.020$ and $-0.018$.
The asterisks indicate that the null hypothesis is
rejected at respectively the 1\% (***), the 5\% (**) or at the 10\% (*) significance level.
Notice that the same decision is reached by all the tests when
$\beta_{20} = -0.018$ but not when
$\beta_{20} = -0.026, -0.024, -0.022$ and $-0.020$. In all cases considered here, the
score and gradient tests lead to the same conclusion. Additionally, the
likelihood ratio and Wald tests display the smallest $p$-values
in all cases, in accordance with their
liberal behaviours observed in our simulation study.
\begin{table}[!htp]
{\footnotesize
\begin{center}
\caption{Test statistics for $\mathcal{H}_0:\beta_{2}=\beta_{20}$ against $\mathcal{H}_1:\beta_{2}\neq\beta_{20}$
($p$-values between parentheses).}\label{tab_new}
\begin{tabular}{cccccc}\hline
          &      \multicolumn{5}{c}{$\beta_{20}$}\\  \cline{2-6}
statistic &$-0.026$        & $-0.024$       &  $-0.022$      & $-0.020$       & $-0.018$      \\ \hline
$S_1$     &$ 7.314\,(0.007)^{***}$ & $5.606\,(0.018)^{**}$ & $4.011\,(0.045)^{**}$ & $2.591\,(0.107)$ & $1.411\,(0.235)$\\
$S_2$     &$11.409\,(0.001)^{***}$ & $8.193\,(0.004)^{***}$ & $5.509\,(0.019)^{**}$ & $3.355\,(0.067)^{*}$ & $1.733\,(0.188)$\\
$S_3$     &$ 5.872\,(0.015)^{**}$ & $4.636\,(0.031)^{**}$ & $3.407\,(0.065)^{*}$ & $2.251\,(0.134)$ & $1.249\,(0.264)$\\
$S_4$     &$ 5.728\,(0.017)^{**}$ & $4.611\,(0.032)^{**}$ & $3.458\,(0.063)^{*}$ & $2.332\,(0.127)$ & $1.321\,(0.250)$\\ \hline
\end{tabular}
\end{center}
}
\end{table}

Notice that the sample size is $n=31$, but if $n$ were smaller,
the tests could lead to different conclusions. To
illustrate this, a randomly chosen subset of the data set with $n=10$ was drawn.
The null hypothesis to be tested is $\mathcal{H}_{0}:\beta_{2}=0$.
The observed value of the test statistics
are $S_1 = 2.939$ ($p$-value: 0.086), $S_2 = 2.980$
($p$-value: 0.084), $S_3 = 2.491$ ($p$-value: 0.114) and
$S_4 = 2.682$ ($p$-value = 0.101).  Hence, at the 10\% significance level, the score and gradient
tests do not reject the null hypothesis unlike the likelihood ratio
and Wald tests, which are much more oversized than the score
and gradient tests as evidenced by our simulation results.

\section{Concluding remarks}\label{conclusions}

The dispersion models (DMs) extend the well-known generalised
linear models \citep{NelderWedderburn1972} and also the exponential family
nonlinear models \citep{CordPaula1989}.
Additionally, the class of DMs covers a comprehensive range of
non-normal distributions. In this paper, we dealt with
the issue of performing hypothesis testing in DMs. We considered the three
classic tests, likelihood ratio, Wald and score tests, and a recently
proposed test, the gradient test. We have derived formulae for the asymptotic
expansions up to order $n^{-1/2}$ of the distribution functions
of the likelihood ratio, Wald, score and gradient
statistics, under a sequence of Pitman alternatives,
for testing a subset of regression  parameters and for
testing the dispersion parameter.
The formulae derived are simple to be used analytically
to obtain closed-form expressions for these expansions in
special models. Also, the power of all four
criteria, which are equivalent to first order,
were compared under specific conditions based on second order
approximations. Additionally, we present Monte Carlo simulations in order to
compare the finite-sample performance of these tests. From the simulation results
we can conclude that the score and gradient tests should be preferred.
Finally, we present an empirical application for illustrative purposes.

\section*{Acknowledgments}

We gratefully acknowledge the financial support of FAPESP and CNPq (Brazil).

\end{document}